\documentclass[12pt,a4j]{article}

\title{A Necessary Condition for existence of Lie Symmetries in
  Quasihomogeneous Systems of Ordinary Differential Equations}
\author{Y.~Hirata$^{\rm a}$\thanks{E-mail: {\tt
      yhirata@allegro.phys.nagoya-u.ac.jp}} \ and K.~Imai$^{\rm
    b}$\thanks{E-mail: {\tt kimai@daido-it.ac.jp}} \\
  $^{\rm a}$Department of Physics, Nagoya University, \\
  Nagoya, 464-8602, Japan \\
  $^{\rm b}$Departments of Engineering, Daido Institute of Technology,
  \\
  Nagoya, 457-8530, Japan}
\date{}

\usepackage{amsmath}
\usepackage{amsthm}

\def\C{\mathbf{C}}
\def\R{\mathbf{R}}
\def\Q{\mathbf{Q}}

\def\N{\mathbf{N}}

\def\del#1{\frac{\partial}{\partial #1}}
\def\pdif#1#2{\frac{\partial #1}{\partial #2}}
\def\round#1{\frac{\partial}{\partial #1}}
\def\Round#1{\partial / \partial #1}

\def\ds{\displaystyle}

\newtheorem{lemma}{Lemma}
\newtheorem{definition}{Definition}
\newtheorem{remark}{Remark}
\newtheorem{theorem}{Theorem}
\newtheorem{proposition}{Proposition}

\def\diag{\mathrm{diag}}

\begin{document}

\maketitle

\begin{abstract}
  Lie symmetries for ordinary differential equations are studied.
  In systems of ordinary differential equations, there do not always
  exist non-trivial Lie symmetries around equilibrium points.
  We present a necessary condition for existence of Lie symmetries
  analytic in the neighbourhood of an equilibrium point.
  In addition, this result can be applied to a necessary condition for
  existence of a Lie symmetry in quasihomogeneous systems of ordinary
  differential equations.
  With the help of our main theorem, it is proved that several systems
  do not possess any analytic Lie symmetries.
\end{abstract}

%%% Local Variables: 
%%% mode: latex
%%% TeX-master: "elsevier"
%%% End: 

%\begin{flushleft}
%PACS codes: 02.20.Sv, 02.30.Hq, 02.90.+p \\
%Keywords: Lie symmetry, Kowalevskaya exponent, quasihomogeneous system
%\end{flushleft}

\section{Introduction}
\label{sec:intro}

A Lie symmetry is defined as a vector field which infinitesimally
transforms a solution of a considered differential equation to, in
general, another solution.
Existence of Lie symmetries of differential equations gives us a lot of
information on dynamics, \textit{e.g.}, existence of similarity
solutions, phase space structure, solvability of the system, and so
on\cite{Olv86,Ibr94}.

For a system of ordinary differential equations, a Lie symmetry can be
regarded as a vector field which commute with the vector field
accompanying with the ordinary differential equations.
However, in general, it is difficult to search Lie symmetries.

Hereby, it is known that existence or non-existence of first integrals
has relevant to Kowalevskaya exponents, first named by Yoshida.
Yoshida argued a necessary condition for existence of first integrals
by using Kowalevskaya exponents\cite{Yos83}.
Later Yoshida's theorem was extended to a necessary condition for
existence of invariant tensor fields with general orders by
Kozlov\cite{Koz92}.
Although above two theorems seem powerful, they have an important weak
point.
Their theorems gives some conditions for Kowalevskaya exponents under
existence of a ``non-degenerate'' first integral (or invariant tensor
field).
Therefore even if the conditions are not satisfied, one can claim only
non-existence of ``non-degenerate'' ones.

Recently Furta\cite{Fur96} and Goriely\cite{Gor96} independently
overcame this weak point although the assertion is a little weak.
Because this method does not suppose non-degeneracy of first
integrals, one can use the method for proofs of non-existence of first
integrals.
In this paper we show a necessary condition for existence of analytic
Lie symmetries by using Furta's method.

This paper is outlined as follows.
Local Lie symmetries always exist around regular points of a system of
ordinary differential equations.
However in the neighbourhood of an equilibrium point, there are not
necessarily non-trivial Lie symmetries.
In the next section we give a necessary condition for existence of Lie
symmetries around an equilibrium point.
We devote Section~\ref{sec:qhsys} to an introduction of
quasihomogeneous (similarity invariant or weighted homogeneous)
systems and Lie symmetries, which are main subjects in this paper.
We give our main theorem and its proof in
Section~\ref{sec:maintheorem}.
Section~\ref{sec:semi} is devoted to extension of the previous section
to semi-quasihomogeneous systems, which is defined in the section.
A concrete example is given in Section~\ref{sec:example}.
We conclude this paper in Section~\ref{sec:conclusion}.

\section{Ordinary differential equations and Lie symmetries}
\label{sec:ODE}

In this section, we present a basic lemma on existence of Lie
symmetries in a system of analytic ordinary differential equations.
It can be regarded as a necessary condition for existence of local Lie 
symmetries around an equilibrium point of a system of analytic
ordinary differential equations.

Assume the system
\begin{equation}
  \label{eq:gsystem}
  \dot{x} = g(x), \quad x \in \C^n
\end{equation}
is analytic in the neighbourhood of the origin $x = 0$.
%Let $f$ be an analytic function $\C^n \to \C^n$ and consider a system
%of ordinary differential equations
%\begin{equation}
%  \label{eq:gsystem}
%  \dot{x} = f(x).
%\end{equation}
We also write the system~\eqref{eq:gsystem} as
\begin{equation}
  \label{eq:gsystemExplicit}
  \dot{x}^i = g^i(x), \quad i = 1, \dots n
\end{equation}
by using the $i$-th component of $x$ and $g$.
Because we are interested in Lie symmetries around an equilibrium
point, we set $g(0) = 0$ without loss of generality.

We define Lie symmetries of the system~\eqref{eq:gsystem} (or
equivalently~\eqref{eq:gsystemExplicit}).
Let $\mathit{X}_g$ be the vector field accompanying with the
system~\eqref{eq:gsystem}, that is,
\[
X_g = \sum_{i=1}^{n} g^i(x) \frac{\partial}{\partial x^i}.
\]
%\begin{equation}
%  \label{eq:fVecFld}
%  X_g = \sum_{i=1}^{n} g^i(x) \frac{\partial}{\partial x^i}.
%\end{equation}
\begin{definition}[Lie symmetries]
  A vector field
  \begin{equation}
    \label{eq:phi}
    X_\phi = \sum_{i=1}^{n} \phi^i(x) \frac{\partial}{\partial x^i}
  \end{equation}
  is called a Lie symmetry of the system~\eqref{eq:gsystem} if the
  vector field~\eqref{eq:phi} commutes with $\mathit{X}_g$ in the
  sense of the Lie bracket, that is,
  \[
  [ X_g, X_\phi ] := X_g X_\phi - X_\phi X_g = 0.
  \]
%  \begin{equation}
%    \label{eq:commute}
%    [ X_g, X_\phi ] := X_g X_\phi - X_\phi X_g = 0.
%  \end{equation}
\end{definition}
\begin{remark}
  In general, one has to set Lie symmetries in the form of
  \[
  X = \tau (x,t) \round{t} + \sum_{i=1}^{n} \phi^i (x,t) \round{x^i}.
  \]
  However, in this paper, we restrict ourselves to autonomous Lie
  symmetries in the form of the expression~\eqref{eq:phi}.
\end{remark}
We say an \textit{analytic} Lie symmetry or a \textit{polynomial} Lie
symmetry if all the components of the vector field $\phi^i(x)$ are
analytic or polynomial, respectively.
In addition, let us define order of Lie symmetries.
\begin{definition}[Order of vector fields]
  An analytic vector field (or Lie symmetry) $X_\phi = \sum \phi^i(x)
  \Round{x^i}$ is of order $k$ if the Taylor expansion of the
  coefficient vector $\phi (x) = {}^t (\phi^1(x), \dots, \phi^n(x))$
  starts from $k$-th order terms.
\end{definition}
%an analytic Lie symmetry $X_\phi$ is called a Lie
%symmetry of degree $m$ if all the components $\phi^i(x)$ start with
%$m$th order terms of $x$.

%Suppose the Jacobian matrix at the origin $\partial g / \partial x (0)
%= A$ is diagonalizable.
%We therefore can set $A = \diag (\lambda_1, \dots, \lambda_n)$.
The following statement holds on an analytic Lie symmetry around an
equilibrium point.
\begin{lemma}
  \label{lem:ode}
  Set $\partial g / \partial x (0) = A$, and $\lambda_1, \dots,
  \lambda_n$ are eigenvalues of $A$.
  Suppose $\det A \ne 0$ and $A$ is diagonalizable.
  If the system~\eqref{eq:gsystem} possesses an analytic Lie symmetry
  of order $k$ in the neighbourhood of the origin, then for some $j
  \in \{ 1, \dots, n\}$, a resonance condition
  \begin{equation}
    \label{eq:resonanceODE}
    \sum_{i=1}^{n} k_i \lambda_i = \lambda_j, \quad k_i \in \N \cup
    \{ 0 \}, \quad \sum_{i=1}^{n} k_i = k
  \end{equation}
  holds.
\end{lemma}
\begin{remark}
  The assumption $\det A \ne 0$ is not necessary.
  In fact, if $\det A = 0$, the set $\lambda_1, \dots, \lambda_n$
  always satisfies at least one resonance relation of the
  type~\eqref{eq:resonanceODE}.
\end{remark}
\begin{remark}
%  On the contrary, the assumption $k \ge 2$ is essential.
  The resonance relation~\eqref{eq:resonanceODE} possesses a trivial
  solution $k_i = \delta_{ij}$.
  This corresponds with the fact that any vector field trivially
  commute with itself.
  Therefore the contrapositive statement of the lemma does not
  tell non-existence of Lie symmetries of order $1$.
%  In fact, if only the trivial resonance exists, Lie symmetries whose
%  linear parts accord with $X_{Ax}$ and nonlinear parts are different
%  from those of $X_g$ can possibly exist.
\end{remark}
\begin{proof}
  Let $X_\phi$ be an analytic Lie symmetry of order $k$.
  Therefore the relation $[ X_g, X_\phi] = 0$ holds.
  $A$ has already been diagonalized and hence $A = \diag (\lambda_1,
  \dots, \lambda_n)$.

  From the analyticity of $X_\phi$, we can expand $X_\phi$ into the
  power series with respect to $x$ as
  \[
  X_\phi = X_{\phi_k} + X_{\phi_{k+1}} + \dots
  \]
%  \begin{equation}
%    \label{eq:expandSymm}
%    X_\phi = X_{\phi_k} + X_{\phi_{k+1}} + \dots
%  \end{equation}
  where
  \[
  X_{\phi_m} = \sum_{j=1}^{n} \phi^j_l (x) \del{x^j}, \quad m = k,
  k+1, \dots
  \]
%  \begin{equation}
%    \label{eq:expandCoef}
%    X_{\phi_l} = \sum_{j=1}^{n} \phi^j_l (x) \del{x^j}, \quad l = k,
%    k+1, \dots
%  \end{equation}
  and $\phi^i_m$ are homogeneous polynomials of order $m$.

%  Consider the first term of the expansion~\eqref{eq:expandSymm}
%  \begin{equation}
%    \label{eq:symmFirstTerm}
%    X_{\phi_{(0)}} = \sum_{i=1}^{n} \phi^i_{(0)} \del{x^i}
%  \end{equation}
%  where $\phi_{(0)} = {}^t (\phi^1_{(0)}, \dots, \phi^n_{(0)})$ is a
%  constant vector.
%  By considering the lowest term of the equation $[X_f, X_\phi] = 0$,
%  we have
%  \begin{equation}
%    \label{eq:LemFirstTerm}
%    [ X_{Ax}, X_{\phi_{(0)}} ] = - \sum_{i,j} A_{i,j} \phi^j_{(0)}
%    \del{x^i} = 0.
%  \end{equation}
%  Hence $A \phi_{(0)} = 0$.
%  Since $\det A \ne 0$, we have $\phi_{(0)} = 0$.

%  Next we assume $X_{\phi_{(0)}} \equiv \dots \equiv
%  X_{\phi_{(k-1)}}$.
%  Hence the lowest term $[ X_{Ax}, X_{\phi_{(k)}}]$ must be zero.
  Since the lowest order terms of $[ X_g, X_\phi ]$ must be zero, we
  obtain the relation
  \[
  [ X_{Ax}, X_{\phi_k} ] = 0.
  \]
%  \begin{equation}
%    \label{eq:lammasymk}
%    [ X_{Ax}, X_{\phi_k} ] = 0.
%  \end{equation}
  Let us explicitly write $\phi^j_k(x)$ as
  \[
  \phi^j_k (x) = \sum_{k_1 + \dots + k_n = k} \phi^j_{k_1 \dots k_n}
  (x^1)^{k_1} \dots (x^n)^{k_n}.
  \]
%  \begin{equation}
%    \label{eq:kthcoef}
%    \phi^j_k (x) = \sum_{k_1 + \dots + k_n = k} \phi^j_{k_1 \dots k_n}
%    (x^1)^{k_1} \dots (x^n)^{k_n}.
%  \end{equation}
  Straightforward computation gives
  \[
  \sum_{j=1}^{n} \sum_{k_1 + \dots + k_n = k}
  \left(
    \sum_{i=1}^{n} \lambda_i k_i - \lambda_j
  \right)
  \phi^j_{k_1 \dots k_n} (x^1)^{k_1} \dots (x^n)^{k_n} \round{x_j} =
  0.
  \]
%  \begin{equation}
%    \label{eq:lemmasymkexp}
%    \sum_{j=1}^{n} \sum_{k_1 + \dots + k_n = k}
%    \left(
%      \sum_{i=1}^{n} \lambda_i k_i - \lambda_j
%    \right)
%    \phi^j_{k_1 \dots k_n} (x^1)^{k_1} \dots (x^n)^{k_n} \round{x_j} =
%    0.
%  \end{equation}
%  \begin{equation}
%    \label{eq:lemmasymkexp}
%    \sum_{j=1}^{n}
%    \left\{
%      \sum_{i=1}^{n} \sum_{k^j_1 + \dots + k^j_n = k} ( \lambda_i
%      k^j_i - \delta_{ij} \lambda_j ) \phi^j_{k^j_1 \dots k^j_n}
%      (x^1)^{k^j_1} \dots (x^n)^{k^j_n}
%    \right\}
%    \round{x_j} = 0.
%  \end{equation}
  Since $X_\phi$ is of order $k$, at least one coefficient
  $\phi^j_{k_1 \dots k_n}$ is not zero.
  This implies the condition~\eqref{eq:resonanceODE}.
\end{proof}

\section{Quasihomogeneous systems and Lie symmetries}
\label{sec:qhsys}

Consider the system
\begin{equation}
  \label{eq:fsystem}
  \dot{x} = f(x), \quad x \in \C^n
\end{equation}
where $f(x)$ is analytic in the neighbourhood of the origin.
We define quasihomogeneous systems, which is also called similarity
invariant systems or weighted homogeneous systems.
\begin{definition}[Quasihomogeneous systems\cite{Yos83}]
  The system~\eqref{eq:fsystem} is called a quasihomogeneous
  system if the system is invariant under the transformation
  \begin{equation}
    \label{eq:trans}
    t \mapsto \alpha^{-1} t, \quad x^i \mapsto \alpha^{g_i} x^i, \quad
    g_i \in \Q, \quad i = 1, \dots, n
  \end{equation}
  for $\forall \alpha \in \R^+$, \textit{i.e.}, $f^i (x),\, i = 1,
  \dots,n$ satisfy
  \begin{equation}
    \label{eq:transf}
    f^i (\alpha^{g_1} x^1, \dots, \alpha^{g_n} x^n) = \alpha^{g_i + 1}
    f^i (x^1, \dots, x^n).
  \end{equation}
\end{definition}
Hereafter we refer $g_1, \dots, g_n$ as the weight exponents.

\begin{remark}
  By defining $G = \diag (g_1, \dots, g_n)$ and $\alpha^G = \diag
  (\alpha^{g_1}, \dots, \alpha^{g_n})$, the expression
  \eqref{eq:transf} can be rewritten as
  \begin{equation}
    \label{eq:transfvec}
    f^i (\alpha^G x) = \alpha^{g_i + 1} f^i (x), \quad i = 1, \dots,
    n.
  \end{equation}
  Moreover, since $f = {}^t (f^1, \dots, f^n)$, we can simply write
  the expression~\eqref{eq:transfvec} as
  \begin{equation}
    \label{eq:transfmat}
    f ( \alpha^G x) = \alpha^{G+E} f(x)
  \end{equation}
  where $E$ is the $n$-dimensional unit matrix.
\end{remark}
%\begin{remark}
%  From the expression~\eqref{eq:transfmat}, if $f$ is analytic in the
%  neighbourhood of the origin, then the origin is an equilibrium
%  point.
%\end{remark}

Choose $l$ be the smallest positive integer such that all the
quantities $lg_1, \dots, lg_n$ are integers.
Then we call the system~\eqref{eq:fsystem} a quasihomogeneous system
of ramification degree $l$.
%Set $\alpha = \mu^l$ and $lG = S$.

Let us define degrees of quasihomogeneous functions and
quasihomogeneous Lie symmetries.
\begin{definition}[Quasihomogeneous functions]
  A function $F(x)$ is called a quasihomogeneous function of degree
  $M$ if the relation
  \begin{equation}
    \label{eq:qhfunc}
    F(\alpha^G x) = \alpha^{M/l} F(x)
%    F(\alpha^G x) = \alpha^{\frac{M}{l}} F(x)
  \end{equation}
  holds.
\end{definition}
\begin{proposition}
  An analytic quasihomogeneous function has an integer degree.
\end{proposition}
\begin{proof}
  Set $S = lG$ and $\alpha^{1/l} = \mu$.
  Hence the expression~\eqref{eq:qhfunc} can be rewritten as
  \begin{equation}
    \label{eq:qhfunctrans}
    F(\mu^S x) = \mu^M F(x).
  \end{equation}
  The l.h.s.\ of the expression~\eqref{eq:qhfunctrans} can be
  expanded into the Laurent series with respect to $\mu$ from the
  analyticity of $F$.
  This implies $M$ must be an integer.
\end{proof}
\begin{definition}[Quasihomogeneous Lie symmetries]
  We call a Lie symmetry $X_\phi$ a quasihomogeneous Lie symmetry of
  degree $M$ if all the coefficients $\phi^i(x), \, i = 1, \dots, n$
  are quasihomogeneous functions of degree $M + g_i l$.
\end{definition}
\begin{remark}
  If $X_\phi$ is a Lie symmetry of degree $M$, the following relation
  holds:
  \[
  X_\phi |_{\alpha^G x} = \alpha^{M/l} X_\phi.
  \]
%  \begin{equation}
%    \label{eq:qhsym}
%    X_\phi |_{\alpha^G x} = \alpha^{M/l} X_\phi.
%  \end{equation}
%  In other words, a quasihomogeneous Lie symmetry of degree $M$ is
%  invariant under the transformation~\eqref{eq:trans} up to the extra
%  coefficient $\alpha^M$.
\end{remark}

For quasihomogeneous systems, it is sufficient to consider
quasihomogeneous Lie symmetries.
In fact, every analytic Lie symmetry $X_\phi$ splits into the power
series of $\alpha$ with quasihomogeneous vector field coefficients by
the transformation~\eqref{eq:trans}, \textit{i.e.},
\[
X_\phi |_{\alpha^G x} = \sum_{M} \alpha^{M/l} X_{\phi_M},
\]
%\begin{equation}
%  \label{eq:symmetrytrans}
%  X_\phi |_{\alpha^G x} = \sum_{M} \alpha^{M/l} X_{\phi_M},
%\end{equation}
where $X_{\phi_M}$ are quasihomogeneous vector fields of degree $M$.
Since $\alpha$ is arbitrary, each of $X_{\phi_M}$ must be a Lie
symmetry.
Hence we restrict ourselves to quasihomogeneous Lie symmetries.

%If the origin is an equilibrium point of the system~\eqref{eq:fsystem}
%(or equivalently a singular point of $X_f$), it must also be a
%singular point of analytic Lie symmetries.
%This implies that $M \ge 0$.

%We introduce important quantities for quasihomogeneous systems.
%First we give a basic fact without a proof.
Quasihomogeneous systems~\eqref{eq:fsystem} possess a particular
solution in the form
\begin{equation}
  \label{eq:SingularSol}
  x(t) = c t^{-G}
\end{equation}
where $c \ne 0$ is a solution vector of the equation
\begin{equation}
  \label{eq:Cdef}
  f(c) + G c = 0.
\end{equation}
%Hereafter we suppose $c \in D \backslash \{ 0 \}$.
Consider variation of the system~\eqref{eq:fsystem} around the
particular solution~\eqref{eq:SingularSol}.
Let us make the following change of variables $x \mapsto u$ such that
\begin{equation}
  \label{eq:transXU}
  x = t^{-G} ( c + u ).
\end{equation}
Substituting the expression~\eqref{eq:transXU} into the
equation~\eqref{eq:fsystem} and using the
relation~\eqref{eq:transfmat}, we have
\begin{equation}
  \label{eq:ftrans1}
  - G t^{-G-E} ( c + u ) + t^{-G} \dot{u} = t^{-G-E} f(c+u).
\end{equation}
Introducing the new independent variable $\tau = \log t$, the
equation~\eqref{eq:ftrans1} can be expressed as an autonomous system
\begin{equation}
  \label{eq:ftrans2}
  u^\prime = Gc + Gu + f(c+u),
\end{equation}
where ${}^\prime = d / d\tau$.
Moreover, from the analyticity of $f$, one can expand $f$ around $x =
c$ and has
\begin{equation}
  \label{eq:fexpand}
  f(c+u) = f(c) + \pdif{f}{x}(c) u + \tilde{f}(u),
\end{equation}
where $\tilde{f}$ stands for higher order terms.
Substituting the equations~\eqref{eq:fexpand} and~\eqref{eq:Cdef} into
the equation~\eqref{eq:ftrans2}, we have
\begin{equation}
  \label{eq:ftrans3}
  u^\prime =
  \left(
    \pdif{f}{x}(c) + G
  \right)
  u + \tilde{f}(u).
\end{equation}
The quantity in the bracket of the r.h.s.\ of the
equation~\eqref{eq:ftrans3} is called the \textit{Kowalevskaya
  matrix}.
\begin{definition}[Kowalevskaya matrix and exponents\cite{Yos83}]
  The $n \times n$ matrix
  \begin{equation}
    \label{eq:Kmat}
    K = \pdif{f}{x}(c) + G
  \end{equation}
  is called the \textit{Kowalevskaya matrix}.
  Moreover the eigenvalues of $K$ is called the \textit{Kowalevskaya
    exponent}s.
\end{definition}
Thus the quasihomogeneous system~\eqref{eq:fsystem} can be expressed
as
\begin{equation}
  \label{eq:Ksystem}
  u^\prime = Ku + \tilde{f}(u)
\end{equation}
around the particular solution~\eqref{eq:transXU}.
In this paper, we assume that Kowalevskaya matrices are
diagonalizable.
\begin{lemma}[\cite{Yos83,Fur96}]
  \label{lem:ke1}
  $\lambda = -1$ is an eigenvalue of a Kowalevskaya matrix for an
  autonomous system.
\end{lemma}
\begin{proof}
  We show that $\lambda = -1$ is an eigenvalue of $K$ defined as the
  expression~\eqref{eq:Kmat} and $q = Gc$ is an eigenvector belonging
  to $\lambda = -1$.
  First by differentiating the equation~\eqref{eq:transfmat} with
  respect to $\alpha$ and substituting $\alpha = 1$ and $x = c$, one
  obtains
  \begin{equation}
    \label{eq:fdifalpha=1}
    \pdif{f}{x}(c) G c = ( G + E ) f(c).
  \end{equation}
  
  On the other hand,
  \begin{equation}
    \label{eq:Kq}
    Kq = KGc =
    \left(
      G + \pdif{f}{x}(c)
    \right)
    Gc.
  \end{equation}
  By substituting the equation~\eqref{eq:fdifalpha=1} into the
  equation~\eqref{eq:Kq}, we have
  \[
  Kq = G^2 c + ( G + E ) f(c) = G ( Gc + f(c) ) + f(c).
  \]
%  \begin{equation}
%    \label{eq:Kq2}
%    Kq = G^2 c + ( G + E ) f(c) = G ( Gc + f(c) ) + f(c).
%  \end{equation}
  Furthermore, since $Gc + f(c) = 0$,
  \[
  Kq = f(c) = -Gc = -q.
  \]
%  \begin{equation}
%    \label{eq:Kq3}
%    Kq = f(c) = -Gc = -q.
%  \end{equation}
  The lemma is proved.
\end{proof}

\section{Main Theorem}
\label{sec:maintheorem}

Now we give our main theorem:
\begin{theorem}[Main Theorem]
  Let the Kowalevskaya matrix $K$ of the quasihomogeneous
  system~\eqref{eq:fsystem} be diagonalizable, and $\lambda_1, \dots,
  \lambda_n$ be the Kowalevskaya exponents.
  Set $\lambda_1 = -1$ and put $\lambda_0 = 1/l$, where $l$ is the
  ramification degree of the system.
  If the quasihomogeneous system has an analytic Lie symmetry of
  degree $M (\ge 0)$, then, for some $j \in \{ 0, 1, \dots, n \}$, the
  resonance condition
  \begin{equation}
    \label{eq:resonanceThm}
    \sum_{i=1}^{n} k_i \lambda_i = l \lambda_j, \quad k_i \in \N \cup
    \{ 0 \}, \quad k_1 \ge M
  \end{equation}
  holds.
\end{theorem}
\begin{proof}
  Let the vector field
  \begin{equation}
    \label{eq:LieSymMain}
    X_\phi = \sum_{i=1}^{n} \phi^i(x) \del{x^i}
  \end{equation}
  be a quasihomogeneous Lie symmetry of degree $M$.
  The quasihomogeneous system~\eqref{eq:fsystem} is transformed into
  the autonomous system~\eqref{eq:Ksystem} by the
  transformation~\eqref{eq:transXU}.
  Applying the transformation~\eqref{eq:transXU} to the Lie
  symmetry~\eqref{eq:LieSymMain}, we have
  \[
  X_\phi = t^{-M/l} \sum_{i=1}^{n} \phi^i(c+u) \round{u^i}.
  \]
%  \begin{equation}
%    \label{eq:LieSymTrans}
%    X_\phi = t^{-M/l} \sum_{i=1}^{n} \phi^i(c+u) \round{u^i}.
%  \end{equation}
  Putting $t^{-1/l} = u^0$, we have
  Thus
  \begin{equation}
    \label{eq:LieSymTrans2}
    X_\phi = (u^0)^M \sum_{i=1}^{n} \phi^i(c+u) \round{u^i},
  \end{equation}
  which is analytic with respect to $(u^0, u^1, \dots, u^n)$.
  Then the new \textit{independent} variable $u^0$ satisfies a linear
  ordinary differential equation
  \[
  (u^0)^\prime = - \frac{1}{l} u^0.
  \]
%  \begin{equation}
%    \label{eq:u0ode}
%    (u^0)^\prime = - \frac{1}{l} u^0.
%  \end{equation}
  A straightforward computation tells that the transformed vector
  field~\eqref{eq:LieSymTrans2} is a Lie symmetry of the extended
  autonomous system
  \[
  \left\{
    \begin{array}{rcl}
      (u^0)^\prime &=& \ds - \frac{1}{l} u^0 \\
      u^\prime &=& Ku + \tilde{f}(u).
    \end{array}
  \right.
  \]
%  \begin{equation}
%    \label{eq:extendsys}
%    \left\{
%      \begin{array}{rcl}
%        (u^0)^\prime &=& \ds - \frac{1}{l} u^0 \\
%        u^\prime &=& Ku + \tilde{f}(u).
%      \end{array}
%    \right.
%  \end{equation}
  Therefore, applying Lemma~\ref{lem:ode}, we have a resonance
  relation
  \begin{equation}
    \label{eq:resthm1}
    -\frac{k_0}{l} + \sum_{i=1}^{n} k_i \lambda_i = \lambda_j, \quad
    k_0, \dots, k_n \in \N \cup \{ 0 \}, \quad k_0 + \sum_{i=1}^{n}
    k_n \ge M,
  \end{equation}
  where $k_0 = M$.
  Multiplying $l$ to the both hands of the first equation in the
  expression~\eqref{eq:resthm1} and rewriting $k_0 + l k_1 \to k_1$
  and $lk_i \to k_i, i = 2, \dots, n$, we obtain the
  expression~\eqref{eq:resonanceThm}.
  This completes the proof.
\end{proof}

\section{Semi-quasihomogeneous systems}
\label{sec:semi}

We reconsider an $n$-dimensional autonomous system
\begin{equation}
  \label{eq:fsystem2}
  \dot{x} = f(x), \quad x \in \C^n
\end{equation}
which are analytic in the neighbourhood of the origin.
Let the $j$-th component of $f$ be expanded into the Maclaurin series
\begin{equation}
  \label{eq:fexpand2}
  f^j(x) = \sum f^j_{k_1 \dots k_n} (x^1)^{k_1} \dots (x^n)^{k_n}.
\end{equation}
\begin{proposition}[\cite{Fur96}]
  If the system~\eqref{eq:fsystem2} is quasihomogeneous with respect
  to the weight exponents $g_1, \dots, g_n \in \Q$, then all the terms
  in the expansion~\eqref{eq:fexpand2} satisfy the equation
  \[
  g_1 k_1 + \dots + g_n k_n = g_j + 1, \quad j = 1, \dots, n.
  \]
%  \begin{equation}
%    \label{eq:qheq}
%    g_1 k_1 + \dots + g_n k_n = g_j + 1, \quad j = 1, \dots, n.
%  \end{equation}
\end{proposition}
Now let us define \textit{semi-quasihomogeneous system}s.
\begin{definition}[Semi-quasihomogeneous systems\cite{Fur96}]
  The system~\eqref{eq:fsystem2} is called a semi-quasihomogeneous
  system if it can be expressed as the form
  \begin{equation}
    \label{eq:sqhsys}
    \dot{x} = f_m(x) + \hat{f}(x),
  \end{equation}
  where $f_m(x)$ defines a quasihomogeneous system of degree $m$ and
  all the terms in the expansion of $\hat{f}^j(x)$, which is the
  $j$-th component of $\hat{f}(x)$, satisfy either
  \begin{equation}
    \label{eq:sqhineqp}
    g_1 k_1 + \dots + g_n k_n > g_j + 1, \quad j = 1, \dots, n
  \end{equation}
  or
  \begin{equation}
    \label{eq:sqhineqn}
    g_1 k_1 + \dots + g_n k_n < g_j + 1, \quad j = 1, \dots, n.
  \end{equation}
  In addition, if the inequality~\eqref{eq:sqhineqp}
  (\textit{resp.}~\eqref{eq:sqhineqn}) holds, we call the
  system~\eqref{eq:sqhsys} a positively (\textit{resp.} negatively)
  semi-quasihomogeneous system.
\end{definition}

Let us give a theorem for existence of Lie symmetries in
semi-\hspace{0cm}quasihomogeneous systems:
\begin{theorem}[Lie symmetries for semi-quasihomogeneous systems]
  \label{thm:semi}
  If a positively (\textit{resp.} negatively) semi-quasihomogeneous
  system~\eqref{eq:sqhsys} possesses an analytic (\textit{resp.}
  polynomial) Lie symmetry, then the truncated system
  \begin{equation}
    \label{eq:truncsys}
    \dot{x} = f_m(x)
  \end{equation}
  has an analytic (\textit{resp.} polynomial) Lie symmetry.
\end{theorem}
\begin{proof}
  Applying the scale transformation~\eqref{eq:trans} to the
  system~\eqref{eq:sqhsys}, we have
  \begin{equation}
    \label{eq:sqhsys1}
    \dot{x} = f_m(x) + \hat{f}(x,\mu),
  \end{equation}
  where $\hat{f}(x,\mu)$ is a formal Taylor series with respect to
  $\mu$ (\textit{resp.} $1/\mu$) without constant terms.

  On the other hand, an analytic (\textit{resp.} polynomial) Lie
  symmetry $X_\phi$ of the system~\eqref{eq:sqhsys} is transformed
  into the form $\mu^M ( X_M + \mu X_{M+1} + \mu^2 X_{M+2} + \dots )$
  (\textit{resp.} $\mu^M ( X_M + \frac{1}{\mu} X_{M-1} +
  \frac{1}{\mu^2} X_{M-2} + \dots + \frac{1}{\mu^{N}} X_{M-N} )$) for
  some $M$, where $X_i$ are quasihomogeneous vector fields of degree
  $i$.
  Hence
  \begin{gather}
    \label{eq:sqhsym1}
    X_M + \mu X_{M+1} + \mu^2 X_{M+2} + \dots \\
    \label{eq:sqhsym2}
    \text{(\textit{resp.} } X_M + \frac{1}{\mu} X_{M-1} +
    \frac{1}{\mu^2} X_{M-2} + \dots + \frac{1}{\mu^{N}} X_{M-N}
    \text{)}
  \end{gather}
  is an analytic (\textit{resp.} polynomial) Lie symmetry of the
  system~\eqref{eq:sqhsys1}.

  The system~\eqref{eq:sqhsys1} approaches the truncated
  system~\eqref{eq:truncsys} as $\mu \to 0$ (\textit{resp.} $\mu \to
  \infty$).
  On the other hand, the Lie symmetry~\eqref{eq:sqhsym1}
  (\textit{resp.}~\eqref{eq:sqhsym2}) approaches $X_M$ simultaneously.
  Thus the truncated system~\eqref{eq:truncsys} possesses a Lie
  symmetry $X_M$.
\end{proof}

\section{Example}
\label{sec:example}

Consider the $2$-dimensional quadratic systems
\begin{subequations}
  \label{eq:2dim2ord}
  \begin{align}
    \dot{x}^1 &= (x^1)^2 + x^1 x^2 \\
    \dot{x}^2 &= a x^1 x^2 + (x^2)^2
  \end{align}
\end{subequations}
with a parameter $a$, which can be regarded as particular cases of
the Lotka--Volterra system.

The weight exponents of the system is $g_1 = g_2 = 1$ and hence the
ramification degree $l = 1$.
%The vector $c$ is easily computed as ${}^t ( c_1, - c_1 - 1 ),\, c_1
%\in \C$ for $a = 1$ and ${}^t (-1, 0), {}^t (0, -1)$ for $a \ne 1$.
The vector $c$ is easily computed as
\[
c = 
\begin{cases}
  {}^t ( c_1, - c_1 - 1 ),\, c_1 \in \C, & a = 1, \\
  {}^t (-1, 0), {}^t (0, -1), & a \ne 1.
\end{cases}
\]
We get the first Kowalevskaya exponent $\lambda_1 = -1$.
The other Kowalevskaya exponent $\lambda_2 = 0$ or $1 - a$.

We first consider the case $\lambda_2 = 0$.
The resonance relation~\eqref{eq:resonanceThm} is written as
\begin{equation}
  \label{eq:res1}
  -k_1 = -1, \quad k_1 \ge 0, \quad k_1 \ge M,
\end{equation}
and
\begin{equation}
  \label{eq:res2}
  -k_1 = 0, \quad k_1 \ge 0, \quad k_1 \ge M,
\end{equation}
respectively.
The equations~\eqref{eq:res1} and~\eqref{eq:res2} implies $k_1 = 0,1$,
and hence $M \le 1$.
In this example $M \ge -1$, because analytic quasihomogeneous Lie
symmetries with constant coefficients have degree $-1$.
Thus we check the case $M = -1, 0, 1$.
Straightforward computation gives the only case $a = 1$ with
non-trivial analytic Lie symmetries
\begin{gather*}
  X_1 = x^1 \round{x^1} - x^1 \round{x^2}, \\
  X_2 = x^2 \round{x^1} - x^2 \round{x^2}, \\
  X_3 = (x^1)^2 \round{x^1} + x^1 x^2 \round{x^2}, \\
  \intertext{and}
  X_4 = x^1 x^2 \round{x^1} + (x^2)^2 \round{x^2}.
\end{gather*}
%\begin{subequations}
%  \label{eq:symmetries}
%  \begin{gather}
%    X_1 = x^1 \round{x^1} - x^1 \round{x^2}, \\
%    X_2 = x^2 \round{x^1} - x^2 \round{x^2}, \\
%    X_3 = (x^1)^2 \round{x^1} + x^1 x^2 \round{x^2}, \\
%    \intertext{and}
%    X_4 = x^1 x^2 \round{x^1} + (x^2)^2 \round{x^2}.
%  \end{gather}
%\end{subequations}
In other words, when $a \ne 1$, there are no analytic Lie symmetries.
Incidentally, although we have the other second Kowalevskaya exponent
$\lambda_2 = 1 - a$ for $a \ne 1$, it is obvious that there exist no
other analytic Lie symmetries.

The $4$-dimensional Lie algebra $\mathcal{L}_4$ spanned by $\{ X_1,
X_2, X_3, X_4 \}$ is solvable and therefore the
system~\eqref{eq:2dim2ord} can be integrated by using $\mathcal{L}_4$.
The commutator table is given in Table~\ref{tab:ComTab}.
\begin{table}[htbp]
  \begin{center}
%      \begin{tabular}{|c||c|c|c|c|}
%        \multicolumn{5}{c}{$X_j$} \\ \hline
%        & $X_1$ & $X_2$ & $X_3$ & $X_4$ \\ \hline\hline
%        $X_1$ & \phantom{$-X_1+$} $0$ & $-X_1-X_2$ & \phantom{$-X_1+$}
%        $X_3$ & \phantom{$-X_1+$} $-X_3$ \\ \hline
%        $X_2$ & \phantom{$-$} $X_1 + X_2$ & $0$ & $X_4$ & $-X_4$ \\ \hline
%        $X_3$ & $-X_3$ & $-X_4$ & $0$ & $0$ \\ \hline
%        $X_4$ & $X_3$ & $X_4$ & $0$ & $0$ \\ \hline
%      \end{tabular}
      \begin{tabular}{r|rrrr}
        & $X_1$ & $X_2$ & $X_3$ & $X_4$ \\ \hline
        $X_1$ & \phantom{$-X_1+$} $0$ & $-X_1-X_2$ & \phantom{$-X_1+$}
        $X_3$ & \phantom{$-X_1+$} $-X_3$ \\
        $X_2$ & \phantom{$-$} $X_1 + X_2$ & $0$ & $X_4$ & $-X_4$ \\
        $X_3$ & $-X_3$ & $-X_4$ & $0$ & $0$ \\
        $X_4$ & $X_3$ & $X_4$ & $0$ & $0$ \\
      \end{tabular}
%      \begin{tabular}{@{\extracolsep{\fill}}|c||c|c|c|c|}
%        \multicolumn{5}{c}{$X_j$} \\ \hline
%        & $X_1$ & $X_2$ & $X_3$ & $X_4$ \\ \hline\hline
%        $X_1$ & $0$ & $-X_1-X_2$ & $X_3$ & $-X_3$ \\ \hline
%        $X_2$ & $X_1 + X_2$ & $0$ & $X_4$ & $-X_4$ \\ \hline
%        $X_3$ & $-X_3$ & $-X_4$ & $0$ & $0$ \\ \hline
%        $X_4$ & $X_3$ & $X_4$ & $0$ & $0$ \\ \hline
%      \end{tabular}
    \caption{The commutator table of $4$-dimensional Lie algebra
        $\mathcal{L}_4$.
      One can take a sequence of sub-algebra $\mathcal{L}_4 \supset
      \mathcal{L}_3 \supset \mathcal{L}_2$, where $\mathcal{L}_3$ and
      $\mathcal{L}_2$ are sub-algebras spanned by $\{ X_2, X_3, X_4
      \}$ (or $\{ X_1, X_3, X_4 \}$) and $\{ X_3, X_4 \}$,
      respectively.
      Hence this Lie algebra is solvable.}
    \label{tab:ComTab}
  \end{center}
\end{table}

Furthermore, by applying Theorem~\ref{thm:semi}, the
$2$-dimensional Lotka--Volterra systems
\begin{align*}
  \dot{x}^1 &= a_{11} x^1 + a_{12} x^2 + (x^1)^2 + x^1 x^2, \\
  \dot{x}^2 &= a_{21} x^1 + a_{22} x^2 + a x^1 x^2 + (x^2)^2
\end{align*}
%\begin{subequations}
%  \label{eq:LV}
%  \begin{align}
%    \dot{x}^1 &= a_{11} x^1 + a_{12} x^2 + (x^1)^2 + x^1 x^2, \\
%    \dot{x}^2 &= a_{21} x^1 + a_{22} x^2 + a x^1 x^2 + (x^2)^2
%  \end{align}
%\end{subequations}
possess no polynomial Lie symmetries except the case $a = 1$.
Moreover the $2$-dimensional replicator systems
\begin{align*}
  \dot{x}^1 &= (x^1)^2 + x^1 x^2 - x^1 ( (x^1)^2 + x^1 x^2 )\\
  \dot{x}^2 &= a x^1 x^2 + (x^2)^2 - x^2 ( a x^1 x^2 + (x^2)^2)
\end{align*}
%\begin{subequations}
%  \label{eq:replicator}
%  \begin{align}
%    \dot{x}^1 &= (x^1)^2 + x^1 x^2 - x^1 ( (x^1)^2 + x^1 x^2 )\\
%    \dot{x}^2 &= a x^1 x^2 + (x^2)^2 - x^2 ( a x^1 x^2 + (x^2)^2)
%  \end{align}
%\end{subequations}
has no analytic Lie symmetries except the case $a = 1$.

\section{Conclusion}
\label{sec:conclusion}

We give a necessary condition for existence of analytic Lie symmetries
in quasihomogeneous systems of ordinary differential equations.
If the resonance equation~\eqref{eq:resonanceThm} possesses only
finite solutions or no solutions, one can argue a question on
non-existence of analytic Lie symmetries in the given systems.
Indeed, we show non-existence of analytic Lie symmetries in the
$2$-dimensional quadratic systems~\eqref{eq:2dim2ord}, which are
particular cases of the Lotka--Volterra system.

The system~\eqref{eq:2dim2ord} with $a = 1$ possesses the solvable Lie
algebra $\mathcal{L}_4$.
On the other hand, the system has a first integral $x_1 / x_2$.
Although our theorem is not lost this ``solvable'' case, Furta's
theorem miss it because the first integral is not analytic in the
neighbourhood of the origin.
Thus our theorem surpasses Furta's theorem under such a situation.

%Incidentally, our main theorem is consistent with Furta's theorem.
%In fact, if there exists an analytic first integral $\Phi$ with a
%non-zero degree, then the
%vector fields $\Phi X_f, \Phi^2 X_f, \dots$ are all analytic Lie
%symmetries.
%Reflecting this fact, if Furta's resonance equation has at least one
%solution, our resonance equation has infinite solutions, and the
%reverse is also true.
%Thus if there exist no first integrals, one can check that the
%considered system possesses analytic Lie symmetries or not.
The resonance equation~\eqref{eq:resonanceThm} has either finite or
infinite solutions.
If the number of the solutions is finite, one can check the considered
system possesses analytic Lie symmetries or not.
On the other hand, if infinite solutions exist, it is difficult to
argue non-existence of analytic Lie symmetries.
The infinity of solutions is caused by solutions of the resonance
equation in Furta's theorem\cite{Fur96}.
In fact, if there exists an analytic first integral $\Phi$ with a
positive degree, then the vector fields $\Phi X_f, \Phi^2 X_f, \dots$
are all analytic Lie symmetries.
Thus our main theorem is consistent with Furta's theorem.

\section*{Acknowledgements}
\label{acknowledge}

The authors would like to thank M.~Ishii and T.~Konishi for fruitful
discussions.

\end{document}